\documentclass[12pt]{article}

\usepackage{amssymb}
\usepackage{amsmath,amsthm}
\usepackage[dvips]{graphicx}
\usepackage{wrapfig}

\newtheorem{thm}{Theorem}[section]
\newtheorem{cor}[thm]{Corollary}
\newtheorem{lem}[thm]{Lemma}
\newtheorem{prop}[thm]{Proposition}
\theoremstyle{definition}
\newtheorem{defn}[thm]{Definition}
\newtheorem{eg}[thm]{Example}
\newtheorem{q}[thm]{Question}
\theoremstyle{remark}

\title{The warping degree of a link diagram}

\author{Ayaka Shimizu \\ \footnotesize{Department of Mathematics, Osaka City University} \\[-3pt]
\footnotesize{Sugimoto, Sumiyoshi-ku Osaka 558-8585, Japan}\\
\footnotesize{Email: ayakaberrycandy@y2.dion.ne.jp}}

\begin{document}

\maketitle

\begin{abstract}
For an oriented link diagram $D$, the warping degree $d(D)$ is the smallest number of crossing changes which are needed to obtain 
a monotone diagram from $D$. 
We show that $d(D)+d(-D)+sr(D)$ is less than or equal to the crossing number of $D$, 
where $-D$ denotes the inverse of $D$ and $sr(D)$ denotes the number of components which have at least one self-crossing. 
Moreover, we give a necessary and sufficient condition for the equality. 
We also consider the minimal $d(D)+d(-D)+sr(D)$ for all diagrams $D$. 
For the warping degree and linking warping degree, 
we show some relations to the linking number, unknotting number, and the splitting number. 
\end{abstract}

\section{Introduction}

The warping degree and a monotone diagram is defined by Kawauchi for an oriented diagram of a knot, 
a link \cite{kawauchi} or a spatial graph \cite{kawauchi2}. 
The warping degree represents such a complexity of a diagram, and depends on the orientation of the diagram. 
For an oriented link diagram $D$, we say that $D$ is monotone if we meet every crossing point as an over-crossing first 
when we travel along all components of the oriented diagram with an order by starting from each base point. 
This notion is earlier used by Hoste \cite{hoste} and by Lickorish-Millett \cite{lickorish-millett} 
in computing polynomial invariants of knots and links. 
The warping degree $d(D)$ of an oriented link diagram $D$ is the smallest number of crossing changes which are needed to obtain 
a monotone diagram from $D$ in the usual way. 
We give the precise definitions of the warping degree and a monotone diagram in Section 2. 
Let $-D$ be the diagram $D$ with orientations reversed for all components, and we call $-D$ the inverse of $D$. 
Let $c(D)$ be the crossing number of $D$. 
We have the following theorem in \cite{shimizu} which is for a knot diagram: 

\phantom{x}
\begin{thm}{\cite{shimizu}}
Let $D$ be an oriented knot diagram which has at least one crossing point. 
Then we have 
$$d(D)+d(-D)+1\leq c(D). $$
Further, the equality holds if and only if $D$ is an alternating diagram. 
\label{dk}
\end{thm}
\phantom{x}

\noindent Let $D$ be a diagram of an $r$-component link ($r\geq 1$). 
Let $D^i$ be a diagram on a knot component $L^i$ of $L$, and we call $D^i$ a component of $D$. 
We define a property of a link diagram as follows: 

\phantom{x}

\begin{defn}
A link diagram $D$ has \textit{property $C$} if every component $D^i$ of $D$ is alternating, and 
the number of over-crossings of $D^i$ is equal to the number of under-crossings of $D^i$ 
in every subdiagram $D^i\cup D^j$ for each $i\neq j$.
\end{defn}
\phantom{x}

\noindent Note that a diagram $D$ has property $C$ if $D$ is an alternating diagram in the case that $r=1$. 
We generalize Theorem \ref{dk} to a link diagram: 

\phantom{x}
\begin{thm}
Let $D$ be an oriented link diagram, and $sr(D)$ the number of components $D^i$ such that $D^i$ has at least one self-crossing.
Then we have 
$$d(D)+d(-D)+sr(D)\leq c(D).$$
Further, the equality holds if and only if $D$ has property $C$.
\label{mainthm}
\end{thm}
\phantom{x}

\begin{figure}[h]
\begin{center}
\includegraphics[width=25mm]{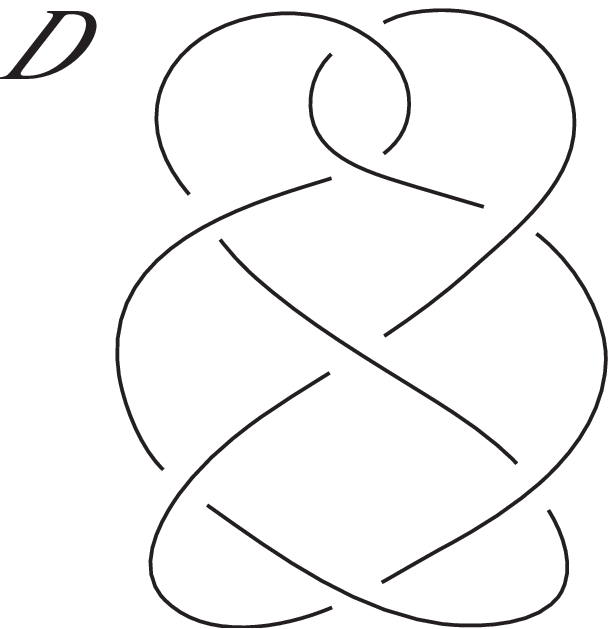}
\end{center}
\caption{}
\label{8-2_5}
\end{figure}

\noindent For example, the link diagram $D$ in Figure \ref{8-2_5} has $d(D)+d(-D)+sr(D)=3+3+2=8=c(D)$.

\noindent Let $D$ be a diagram of a link. 
Let $u(D)$ be the unlinking number of $D$. 
As a lower bound for the value $d(D)+d(-D)+sr(D)$, we have the following inequality: 

\phantom{x}
\begin{thm}
We have 
$$2u(D)+sr(D)\leq d(D)+d(-D)+sr(D).$$
\end{thm}
\phantom{x}

\noindent The rest of this paper is organized as follows. 
In Section 2, we define the warping degree $d(D)$ of an oriented link diagram $D$. 
In Section 3, we define the linking warping degree $ld(D)$, and consider the value $d(D)+d(-D)$ to prove Theorem \ref{mainthm}. 
In Section 4, we show relations of the linking warping degree and the linking number. 
In Section 5, we apply the warping degree to a link itself. 
In Section 6, we study relations to unknotting number and crossing number. 
In Section 7, we define the splitting number and consider relations between the warping degree and the splitting number. 
In Section 8, we show methods for calculating the warping degree and the linking warping degree.

\section{The warping degree of an oriented link diagram}

Let $L$ be an $r$-component link, and $D$ a diagram of $L$. 
We take a sequence ${\bf a}$ of base points $a_i$ ($i=1,2,\dots ,r$), where every component has just one base point 
except at crossing points. 
Then $D_{\bf a}$, the pair of $D$ and ${\bf a}$, is represented 
by $D_{\bf a}=D^1_{a_1}\cup D^2_{a_2}\cup \dots \cup D^r_{a_r}$ 
with the order of ${\bf a}$. 
A self-crossing point $p$ of $D^i_{a_i}$ is a \textit{warping crossing point of} $D^i_{a_i}$ 
if we meet the point first at the under-crossing when we go along the oriented diagram $D^i_{a_i}$ by starting from $a_i$ ($i=1,2,\dots ,r$). 
A crossing point $p$ of $D^i_{a_i}$ and $D^j_{a_j}$ is a \textit{warping crossing point between} $D^i_{a_i}$ \textit{and} $D^j_{a_j}$ 
if $p$ is the under-crossing of $D^i_{a_i}$ ($1\leq i<j\leq r$). 
A crossing point $p$ of $D_{\bf a}$ is a \textit{warping crossing point of} $D_{\bf a}$ 
if $p$ is a warping crossing point of $D^i_{a_i}$ or a warping crossing point between $D^i_{a_i}$ and $D^j_{a_j}$ \cite{kawauchi}.

\begin{figure}[h]
\begin{center}
\includegraphics[width=35mm]{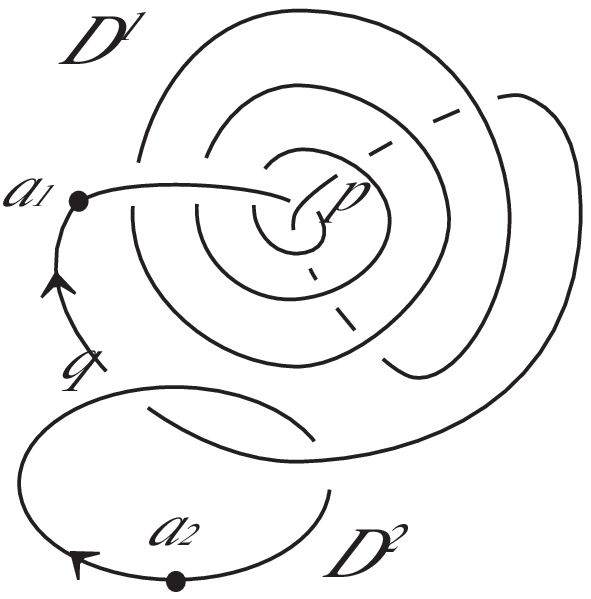}
\end{center}
\caption{}
\label{wcp}
\end{figure}

\noindent For example in Figure \ref{wcp}, $p$ is a warping crossing point of $D^1_{a_1}$, 
and $q$ is a warping crossing point between $D^1_{a_1}$ and $D^2_{a_2}$. 
We define the warping degree for an oriented link diagram \cite{kawauchi}. 
The \textit{warping degree} of $D_{\bf a}$, denoted by $d(D_{\bf a})$, 
is the number of warping crossing points of $D_{\bf a}$. 
The \textit{warping degree} of $D$, denoted by $d(D)$, is the minimal warping degree $d(D_{\bf a})$ 
for all base point sequences ${\bf a}$ of $D$. 
Ozawa showed that a non-trivial link which has a diagram $D$ with $d(D)=1$ is 
a split union of a twist knot or the Hopf link and $r$ trivial knots ($r\geq 0$) \cite{ozawa}. 
Fung also showed that a non-trivial knot which has a diagram $D$ with $d(D)=1$ is a twist knot \cite{stoimenow}.

\phantom{x}

\noindent For an oriented link diagram and its base point sequence $D_{\bf a}=D^1_{a_1}\cup D^2_{a_2}\cup \dots \cup D^r_{a_r}$, 
we denote by $d(D^i_{a_i})$ the number of warping crossing points of $D^i_{a_i}$. 
We denote by $d(D^i_{a_i},D^j_{a_j})$ the number of warping crossing points between $D^i_{a_i}$ and $D^j_{a_j}$. 
By definition, we have that 
\begin{align*}
d(D_{\bf a})=\sum _{i=1}^r d(D^i_{a_i})+\sum _{i<j} d(D^i_{a_i},D^j_{a_j}).
\end{align*}
\noindent Thus, the set of the warping crossing points of $D_{\bf a}$ is divided into two types in the sense 
that the warping crossing point is self-crossing or not. 

\begin{figure}[h]
\begin{center}
\includegraphics[width=35mm]{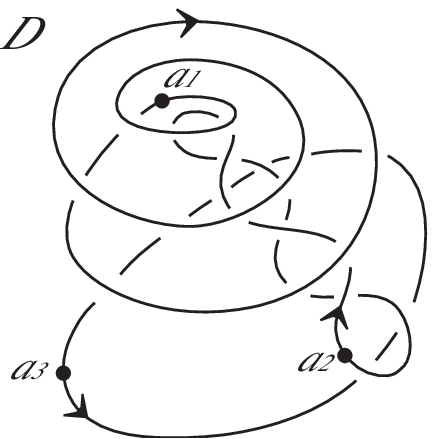}
\end{center}
\caption{}
\label{monotone}
\end{figure}

\noindent The pair $D_{\bf a}$ is \textit{monotone} if $d(D_{\bf a})=0$. 
For example, $D_{\bf a}$ depicted in Figure \ref{monotone} is monotone. 

\noindent Note that a monotone diagram is a diagram of a trivial link. 
Hence we have $u(D)\leq d(D)$, where $u(D)$ is the unlinking number of $D$ (\cite{nakanishi}, \cite{taniyama}).

\section{Proof of Theorem \ref{mainthm}}

\begin{figure}
\begin{center}
\includegraphics[width=70mm]{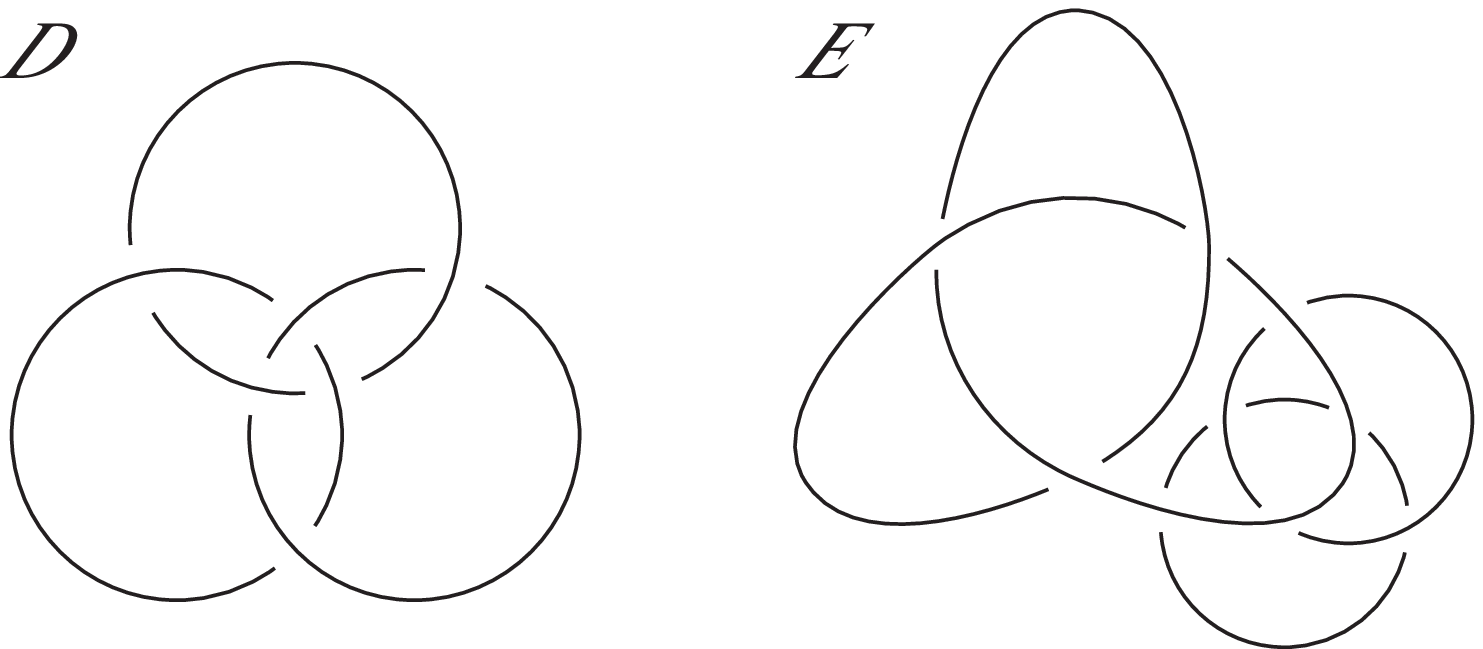}
\end{center}
\caption{}
\label{stacked}
\end{figure}

In this section, we prove Theorem \ref{mainthm}. 
We first define the linking warping degree, which is like a restricted warping degree 
and which has relations to the crossing number and the linking number (see also Section 4). 
The number of non-self warping crossing points does not depend on the orientation. 
We define the \textit{linking warping degree} of $D_{\bf a}$, denoted by $ld(D_{\bf a})$, by the following formula: 

$$ld(D_{\bf a})=\sum _{i<j}d(D^i_{a_i},D^j_{a_j})=d(D_{\bf a})-\sum _{i=1}^rd(D_{a_i}^i),$$
where $D^i_{a_i}, D^j_{a_j}$ are components of $D_{\bf a}$ $(1\leq i<j\leq r)$. 
The \textit{linking warping degree} of $D$, 
denoted by $ld(D)$, is the minimal $ld(D_{\bf a})$ for all base point sequences ${\bf a}$. 
It does not depend on any choices of orientations of components. 
For example, the diagram $D$ in Figure \ref{stacked} has $ld(D)=2$. 
A pair $D_{\bf a}$ is \textit{stacked} if $ld(D_{\bf a})=0$. 
A diagram $D$ is \textit{stacked} if $ld(D)=0$. 
For example, the diagram $E$ in Figure \ref{stacked} is a stacked diagram. 
We remark that a similar notion is mentioned in \cite{hoste}. 
Note that a monotone diagram is a stacked diagram. 
A link $L$ is \textit{completely splittable} if $L$ has a diagram $D$ without non-self crossings. 
Notice that a completely splittable link has some stacked diagrams. 

\phantom{x}

\noindent The \textit{linking crossing number} of $D$, denoted by $lc(D)$, is the number of non-self crossing points of $D$. 
Remark that $lc(D)$ is always even. 
For an unordered diagram $D$, we assume that $D^i$ and $D^i\cup D^j$ denote subdiagrams of $D$ with an order. 
We have the following relation of linking warping degree and linking crossing number.

\phantom{x}
\begin{lem}
We have 
$$ld(D)\leq \frac{lc(D)}{2}.$$
Further, the equality holds if and only if 
the number of over-crossings of $D^i$ is equal to the number of under-crossings of $D^i$ 
in every subdiagram $D^i\cup D^j$ for every $i\neq j$. 

\begin{proof}
Let ${\bf a}$ be a base point sequence of $D$, 
and $\tilde{{\bf a}}$ the base point sequence ${\bf a}$ with the order reversed. 
We call $\tilde{{\bf a}}$ the reverse of ${\bf a}$. 
Since we have that $ld(D_{\bf a})+ld(D_{\tilde{{\bf a}}})=lc(D)$, 
we have the inequality $ld(D)\leq lc(D)/2$. 
Let $D$ be a link diagram such that the number of over-crossings of $D^i$ is equal to 
the number of under-crossings of $D^i$ in every subdiagram $D^i\cup D^j$ for each $i\neq j$. 
Then we have $ld(D_{\bf a})=lc(D)/2$ for every base point sequence ${\bf a}$. 
Hence we have $ld(D)=lc(D)/2$. 
On the other hand, we consider the case the equality $2ld(D)=lc(D)$ holds. 
For an arbitrary base point sequence ${\bf a}$ of $D$ and its reverse $\tilde{{\bf a}}$, 
we have 
$$ld(D_{\bf a})\geq ld(D)=lc(D)-ld(D)\geq lc(D)-ld(D_{\bf a})=ld(D_{\tilde{{\bf a}}})\geq ld(D).$$
Then we have $lc(D)-ld(D_{\bf a})=ld(D)$. 
Hence we have $ld(D_{\bf a})=ld(D)$ for every base point sequence ${\bf a}$. 
Let ${\bf a}'=(a_1,a_2,\dots ,a_{k-1},a_{k+1},a_k,a_{k+2},\dots ,a_r)$ be the base point sequence 
which is obtained from ${\bf a}=(a_1,a_2,\dots ,a_k,a_{k+1},\dots ,a_r)$ 
by exchanging $a_k$ and $a_{k+1}$ ($k=1,2,\dots ,r-1$). 
Then, the number of over-crossings of $D^k$ is equal to the number of under-crossings 
of $D^k$ in the subdiagram $D^k\cup D^{k+1}$ of $D_{\bf a}$ because we have $ld(D_{\bf a})=ld(D_{{\bf a}'})$. 
This completes the proof. 
\end{proof}

\label{2dc}
\end{lem}
\phantom{x}

\noindent We next consider the value $d(D)+d(-D)$ for an oriented link diagram $D$ and the inverse $-D$. 
We have the following proposition: 

\phantom{x}
\begin{prop}
Let $D$ be an oriented link diagram. 
The value $d(D)+d(-D)$ does not depend on the orientation of $D$. 
\begin{proof}
Let $D'$ be $D$ with the same order and another orientation. 
Since we have $d({D^i}')=d(D^i)$ or $d({D^i}')=d(-D^i)$, 
we have $d({D^i}')+d(-{D^i}')=d(D^i)+d(-D^i)$ for each $D^i$ and ${D^i}'$. 
Then we have 

\begin{align*}
d(D')+d(-D')&=\sum _{i=1}^{r}d({D^i}')+ld(D')+\sum _{i=1}^{r}d(-{D^i}')+ld(-D')\\
&=\sum _{i=1}^{r}\{ d({D^i}')+d(-{D^i}')\} +2ld(D')\\
&=\sum _{i=1}^{r}\{ d(D^i)+d(-D^i)\} +2ld(D)\\
&=\sum _{i=1}^{r}d(D^i)+ld(D)+\sum _{i=1}^{r}d(-D^i)+ld(-D)\\
&=d(D)+d(-D).
\end{align*}

\end{proof}
\end{prop}
\phantom{x}

\noindent A link diagram is a \textit{self-crossing diagram} if every component of $D$ has at least one self-crossing. 
In other words, a diagram $D$ of an $r$-component link $L$ is a self-crossing diagram if $sr(D)=r$. 
We have the following lemma: 

\phantom{x}
\begin{lem}
Let $D$ be a self-crossing diagram of an $r$-component link. 
Then we have 
$$d(D)+d(-D)+r\leq c(D).$$
Further, the equality holds if and only if $D$ has property $C$. 

\phantom{x}
\begin{proof}
We have 
\begin{align*}
d(D)+d(-D)+r&=\sum _{i=1}^{r}d(D^i)+ld(D)+\sum _{i=1}^{r}d(-D^i)+ld(-D)+r\\
&=\sum _{i=1}^r\{ d(D^i)+d(-D^i)+1\} +2ld(D)\\
&\leq \sum _{i=1}^r c(D^i)+2ld(D)\\
&\leq \sum _{i=1}^r c(D^i)+lc(D)\\
&=c(D),
\end{align*}
where the first inequality is obtained by Theorem \ref{dk}, 
and the second inequality is obtained by Lemma \ref{2dc}. 
Hence we have the inequality. 
The equality holds if and only if $D$ has property $C$ which is obtained by Theorem \ref{dk} and Lemma \ref{2dc}. 

\end{proof}
\label{self-crossing}
\end{lem}
\phantom{x}

\noindent We give an example of Lemma \ref{self-crossing}.

\phantom{x}
\begin{eg}
In Figure \ref{parallel}, there are three diagrams with 12 crossings. 
The diagram $D$ is a diagram such that any component is alternating 
and has 3 over-non-self crossings and 3 under-non-self crossings. 
Then we have $d(D)+d(-D)+r=12=c(D)$. 
The diagram $D'$ is a diagram which has a non-alternating component diagram. 
Then we have $d(D')+d(-D')+r=10<c(D')$. 
The diagram $D''$ is a diagram such that a component has 2 over-non-self crossings and 4 under-non-self crossings. 
Then we have $d(D'')+d(-D'')+r=10<c(D'')$. 

\begin{figure}[h]
\begin{center}
\includegraphics[width=100mm]{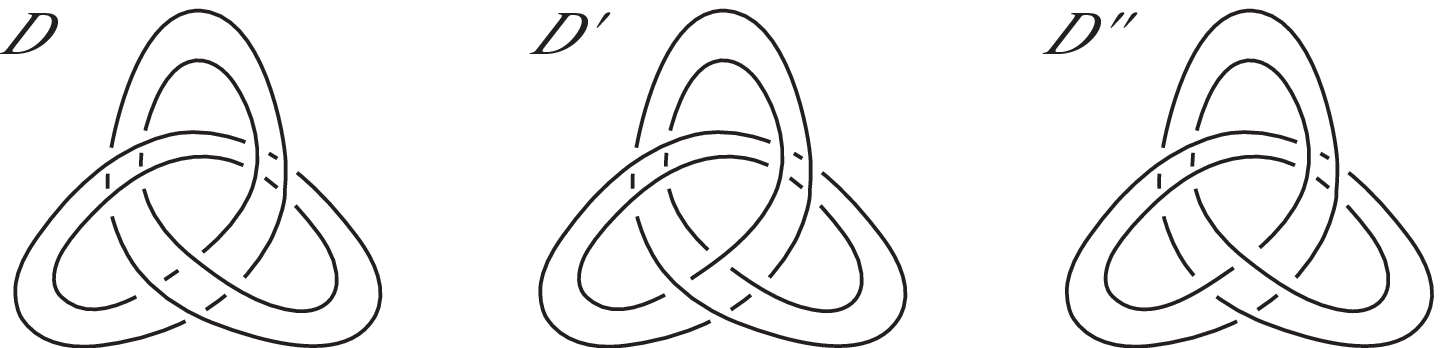}
\end{center}
\caption{}
\label{parallel}
\end{figure}
\end{eg}
\phantom{x}

\noindent Lemma \ref{self-crossing} is only for self-crossing link diagrams. 
We prove Theorem \ref{mainthm} which is for every link diagram.  

\phantom{x}

\begin{figure}[h]
\begin{center}
\includegraphics[width=80mm]{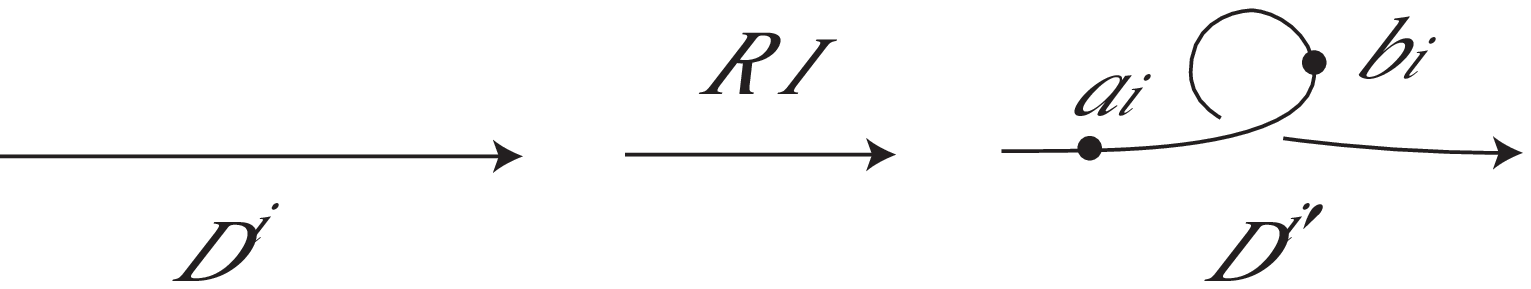}
\end{center}
\caption{}
\label{lr-1}
\end{figure}

\noindent {\it Proof of Theorem \ref{mainthm}.}
For every component $D^i$ such that $D^i$ has no self-crossings, we apply a Reidemeister move of type I as shown in 
Figure \ref{lr-1}. 
Then we obtain the diagram ${D^i}'$ from $D^i$, and ${D^i}'$ satisfies 
$d({D^i}')=d({-D^i}')=0=d(D^i)=d(-D^i)$ and $c({D^i}')=1=c(D^i)+1$. 
For example the base points $a_i$, $b_i$ in Figure \ref{lr-1} satisfy $d(D_{a_i}^i)=d(D^i)=0$, $d(-D_{b_i}^i)=d(-D^i)=0$. 
We remark that every $D^i$ and ${D^i}'$ are alternating. 
We denote by $D'$ the diagram obtained from $D$ by this procedure. 
Since every component has at least one self-crossing, we apply Lemma \ref{self-crossing} to $D'$. Then we have  
$$d(D')+d(-D')+r\leq c(D').$$
And we obtain 
$$d(D)+d(-D)+r\leq c(D)+(r-sr(D)).$$
Hence we have 
$$d(D)+d(-D)+sr(D)\leq c(D).$$
The equality holds if and only if $D$ has property $C$. 
\hfill$\square$

\section{The linking warping degree and linking number}

In this section, we consider the relation of the linking warping degree and the linking number. 
For a crossing point $p$ of an oriented diagram, $\varepsilon (p)$ denotes the sign of $p$, 
namely $\varepsilon (p)=+1$ if $p$ is a positive crossing, and $\varepsilon (p)=-1$ if $p$ is a negative crossing. 
For an oriented subdiagram $D^i\cup D^j$, the \textit{linking number} of $D^i$ with $D^j$ is defined to be 
$$\mathrm{Link}(D^i,D^j)=\frac{1}{2}\sum _{p\in D^i\cap D^j}\varepsilon (p).$$
The linking number of $D^i$ with $D^j$ is independent of the diagram (cf. \cite{cromwell}, \cite{kawauchi}).
We have a relation of the linking warping degree and the linking number of a link diagram in the following proposition: 

\phantom{x}
\begin{prop}
For a link diagram $D$, we have the following (i) and (i\hspace{-1pt}i). 

\begin{description}
\item[(i)] We have 
$$\sum _{i<j} |\mathrm{Link}(D^i,D^j)|\leq ld(D).$$
Further, the equality holds if and only if under-crossings of $D^i$ in $D^i\cup D^j$ are 
all positive or all negative with an orientation for every subdiagram $D^i\cup D^j$ ($i<j$). 

\item[(i\hspace{-1pt}i)] We have 
\begin{align}
\sum _{i<j} |\mathrm{Link}(D^i,D^j)|\equiv ld(D) \ (\bmod \ 2).
\label{mod2}
\end{align}
\end{description}

\phantom{x}
\begin{proof}
\begin{description}
\item[(i)] For a subdiagram $D^i\cup D^j$ ($i<j$) with $d(D^i,D^j)=m$, we show that  
$$|\mathrm{Link}(D^i,D^j)|\leq d(D^i,D^j).$$
Let $p_1, p_2, \dots ,p_m$ be the warping crossing points between $D^i$ and $D^j$, and 
$\varepsilon (p_1)$,$\varepsilon (p_2)$,$\dots ,\varepsilon (p_m)$ the signs of them. 
Since a stacked diagram is a diagram of a completely splittable link, we have 
\begin{align}
\mathrm{Link}(D^i,D^j)-(\varepsilon (p_1)+\varepsilon (p_2)+\dots +\varepsilon (p_m))=0
\label{sign-o-u}
\end{align}
by applying crossing changes at $p_1, p_2, \dots ,p_m$ for $D^i\cup D^j$. 
Then we have 
$$|\mathrm{Link}(D^i,D^j)|=|\varepsilon (p_1)+\varepsilon (p_2)+\dots +\varepsilon (p_m)|\leq m=d(D^i,D^j).$$
Hence we obtain
$$\sum _{i<j} |\mathrm{Link}(D^i,D^j)|\leq ld(D).$$
The equality holds if and only if under-crossings of $D^i$ in $D^i\cup D^j$ are 
all positive or all negative with an orientation for every subdiagram $D^i\cup D^j$ ($i<j$). \\

\item[(i\hspace{-1pt}i)] By the above equality (\ref{sign-o-u}), 
we observe that 
$\mathrm{Link}(D^i,D^j)=\varepsilon (p_1)+\varepsilon (p_2)+\dots +\varepsilon (p_m)=\varepsilon (q_1)+\varepsilon (q_2)+\dots +\varepsilon (q_n)$, 
where $p_k$ (resp. $q_k$) is an under-crossing (resp. an over-crossing) of $D^i$ in $D^i\cup D^j$, $ld(D^i\cup D^j)=m$ and 
$lc(D^i\cup D^j)=m+n$. 
A similar fact is also mentioned in \cite{rolfsen}. We have 
\begin{align*}
\mathrm{Link}(D^i,D^j)&=\varepsilon (p_1)+\varepsilon (p_2)+\dots +\varepsilon (p_m)\\
&\equiv m \ (\bmod \ 2)\\
&=d(D^i,D^j). 
\end{align*}
Hence we have the modular equality 
\begin{align*}
\sum _{i<j} |\mathrm{Link}(D^i,D^j)|\equiv ld(D) \ (\bmod \ 2).
\end{align*}
\end{description}
\end{proof}
\label{lld}
\end{prop}
\phantom{x}

\begin{eg}

In Figure \ref{ex3components}, $D$ has $(0,2,3)$, $E$ has $(0,2,2)$, and $F$ has $(4,4,4)$, 
where $(l,m,n)$ of $D$ denotes that $\sum_{i<j}|\mathrm{Link}(D^i,D^j)|=l$, $ld(D)=m$, and $lc(D)/2=n$.

\begin{figure}[h]
\begin{center}
\includegraphics[width=100mm]{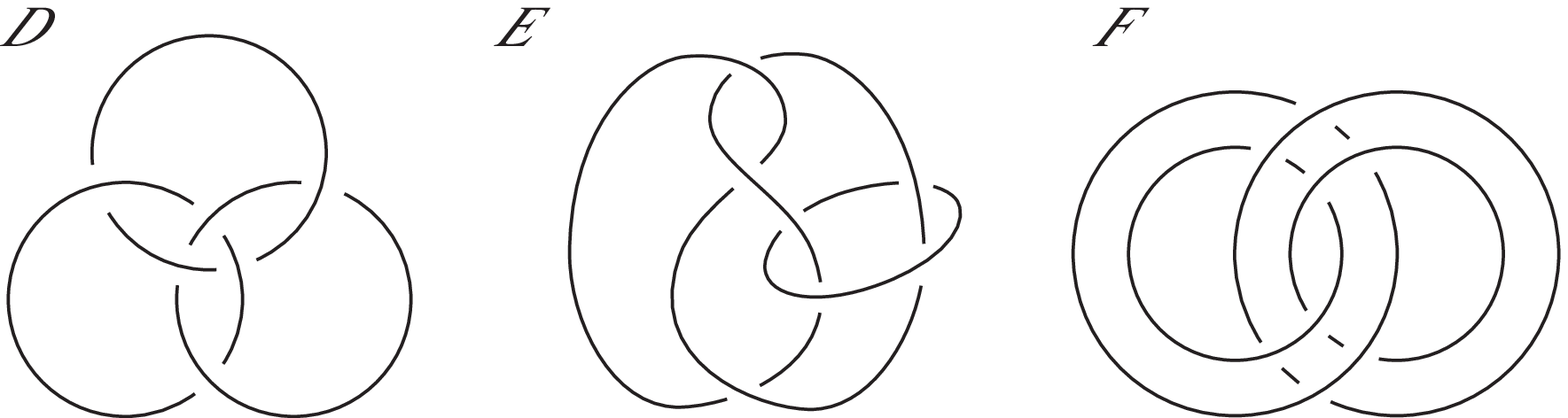}
\end{center}
\caption{}
\label{ex3components}
\end{figure}
\end{eg}

\noindent The \textit{total linking number} of an oriented link $L$ is defined to be $\sum _{i<j}\mathrm{Link}(D^i,D^j)$ with a diagram and an order. 
We have the following corollary: 

\phantom{x}
\begin{cor}
We have 
$$\sum _{i<j}\mathrm{Link}(D^i,D^j)=\sum _{k=1}^r \{ \varepsilon (p_k) | p_k : \textit{a non-self warping crossing point of }D_{\bf a}\},$$
where ${\bf a}$ is a base point sequence of $D$. 
\label{total}
\end{cor}
\phantom{x}

\begin{figure}[h]
\begin{center}
\includegraphics[width=35mm]{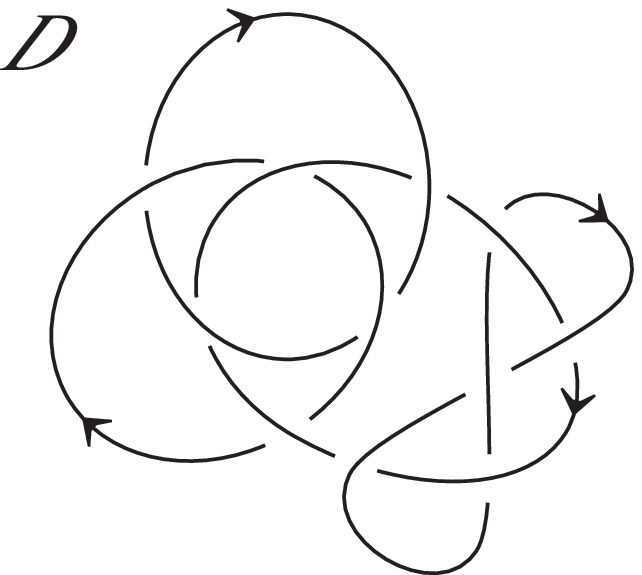}
\end{center}
\caption{}
\label{linkcor}
\end{figure}

\noindent Corollary \ref{total} is useful in calculating the total linking number of a diagram. 
For example in Figure \ref{linkcor}, the diagram $D$ with $4$ components and $11$ crossing points has $ld(D)=4$. 
We have that the total linking number of $D$ is $0$ by summing the signs of only $4$ crossing points. 

\section{To a link invariant}

In this section, we consider the minimal $d(D)+d(-D)$ for minimal crossing diagrams $D$ of $L$ in the following formula: 
$$e(L)=\min \{ d(D)+d(-D) | D:\text{ a diagram of }L\text{ with }c(D)=c(L)\} ,$$
where $c(L)$ denotes the crossing number of $L$. 
In the case where $K$ is a non-trivial knot, we have 
\begin{align}
e(K)+1\leq c(K).
\label{ek}
\end{align}
Further, the equality holds if and only if $K$ is a prime alternating knot \cite{shimizu}. 
Note that the condition for the equality of (\ref{ek}) requires that $D$ is a minimal crossing diagram in the definition of $e(L)$. 
We next define $c^*(L)$ and $e^*(L)$ as follows: 
$$c^*(L)=\min \{c(D) | D:\text{ a self-crossing diagram of }L \} ,$$
$$e^*(L)=\min \{d(D)+d(-D) | D:\text{ a self-crossing diagram of }L\text{ with }c(D)=c^*(L)\} .$$
As a generalization of the above inequality (\ref{ek}), we have the following theorem: 

\phantom{x}
\begin{thm}
For an $r$-component link $L$, we have 
$$e^*(L)+r\leq c^*(L).$$
Further, the equality holds if and only if every self-crossing diagram $D$ of $L$ with $c(D)=c^*(L)$ has property $C$. 

\begin{proof}
Let $D$ be a self-crossing diagram of $L$ with $c(D)=c^*(D)$. 
We assume that $D$ satisfies the equality $d(D)+d(-D)=e^*(L)$. 
Then we have 

\begin{align*}
e^*(L)+r&=d(D)+d(-D)+r\\
&=\sum _{i=1}^{r}d(D^i)+ld(D)+\sum _{i=1}^{r}d(-D^i)+ld(-D)+r\\
&=\sum _{i=1}^r\{ d(D^i)+d(-D^i)+1\} +2ld(D)\\
&\leq \sum _{i=1}^r c(D^i)+2ld(D)\\
&\leq \sum _{i=1}^r c(D^i)+lc(D)\\
&=c(D)=c^*(L),
\end{align*}
where the first inequality is obtained by Theorem \ref{dk}, 
and the second inequality is obtained by Lemma \ref{2dc}. 
If $D$ has a non-alternating component $D^i$, or $D$ has a diagram $D^i\cup D^j$ such that 
the number of over-crossings of $D^i$ is not equal to the number of under-crossings of $D^i$, then 
we have $e^*(L)+r<c^*(L)$. 
On the other hand, the equality holds if $D$ has property $C$. 
\end{proof}
\end{thm}
\phantom{x}

\noindent We have the following example: 

\phantom{x}
\begin{eg}
For non-trivial prime alternating knots $L^1,L^2,\dots ,L^r$ ($r\geq 2$), we have a non-splittable link $L$ by performing 
$n_i$-full twists for every $L^i$ and $L^{i+1}$ ($i=1,2,\dots ,r$) with $L^{r+1}=L^1$ as shown in Figure \ref{full}, 
where we assume that $n_1$ and $n_r$ have the same sign. 

\begin{figure}[h]
\begin{center}
\includegraphics[width=80mm]{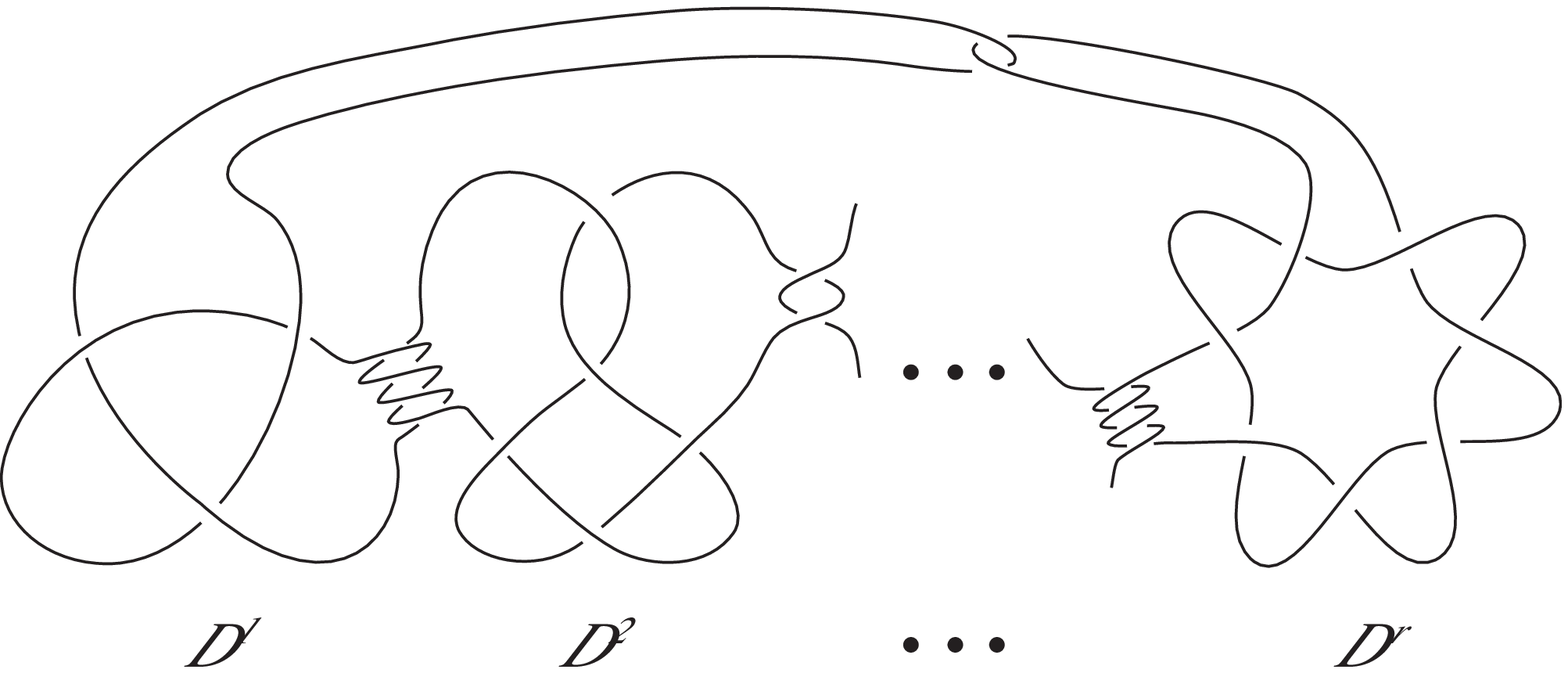}
\end{center}
\caption{}
\label{full}
\end{figure}

\noindent Note that we do not change the type of knot components $L^i$. 
Let $D$ be a diagram of $L$ with $c(D)=c(L)$. 
Then we notice that $D$ is a self-crossing diagram with $c(D)=c^*(L)$. 
We also notice that $D$ has property $C$ 
because $lc(D^i\cup D^j)=2|n_i|$ and $\mathrm{Link}(D^i,D^j)=n_i$, 
and $lc(D^1\cup D^r)=2|n_1+n_r|$ and $\mathrm{Link}(D^1,D^r)=n_1+n_r$ in the case where $r=2$. 
Hence we have $e^*(L)+r=c^*(L)$ in this case. 
\end{eg}
\phantom{x}

\noindent We have the following corollary: 

\phantom{x}
\begin{cor}
Let $L$ be an $r$-component link whose all components are non-trivial. 
Then we have 
$$e(L)+r\leq c(L).$$
Further, the equality holds if and only if every diagram $D$ of $L$ with $c(D)=c(L)$ has property $C$. 

\begin{proof}
Since every diagram $D$ of $L$ is a self-crossing diagram, we have $e(L)=e^*(L)$ and $c(L)=c^*(L)$. 
\end{proof}
\end{cor}
\phantom{x}

\noindent We also consider the minimal $d(D)+d(-D)+sr(D)$ and the minimal $sr(D)$ for diagrams $D$ of $L$ in the following formulas: 
$$f(L)=\min \{ d(D)+d(-D)+sr(D) | D:\text{ a diagram of }L\},$$
$$sr(L)=\min \{ sr(D) | D:\text{ a diagram of }L\}.$$

\noindent Note that the value $f(L)$ and $sr(L)$ also do not depend on the orientation of $L$. 
Jin and Lee mentioned in \cite{jinlee} that every link has a diagram which restricts to a minimal crossing diagram 
for each component. 
Then we have the following proposition: 

\phantom{x}
\begin{prop}
The value $sr(L)$ is equal to the number of non-trivial knot components of $L$. 
\end{prop}
\phantom{x}

\noindent The following corollary is directly obtained from Theorem \ref{mainthm}. 

\phantom{x}
\begin{cor}
We have 
$$f(L)\leq c(L).$$
\begin{proof}
For a diagram $D$ with $c(D)=c(L)$, we have 
$$f(L)\leq d(D)+d(-D)+sr(D)\leq c(D)=c(L),$$
where the second inequality is obtained by Theorem \ref{mainthm}. 
\end{proof}
\label{fcl}
\end{cor}
\phantom{x}

\noindent We have the following question: 

\phantom{x}
\begin{q}
When does the equality $f(L)=c(L)$ hold?
\end{q}
\phantom{x}

\begin{eg}
In Figure \ref{ex-fc}, there are two link diagrams $D$ and $E$. 
We assume that $D$ (resp. $E$) is a diagram of a link $L$ (resp. $M$). 
We have $f(L)=c(L)=5$ because we have $d(D)+d(-D)+sr(D)=2+2+1$ and we know $d(D^i)\geq u(3_1)=1$, $ld(D)\geq 1$, and $sr(D)\geq sr(L)=1$, 
where $D^i$ is any diagram of $3_1$. 
On the other hand, we have that $f(M)<c(M)$ because $f(M)\leq d(E)+d(-E)+sr(E)=3+3+1=7<10=c(M)$. 

\begin{figure}[h]
\begin{center}
\includegraphics[width=70mm]{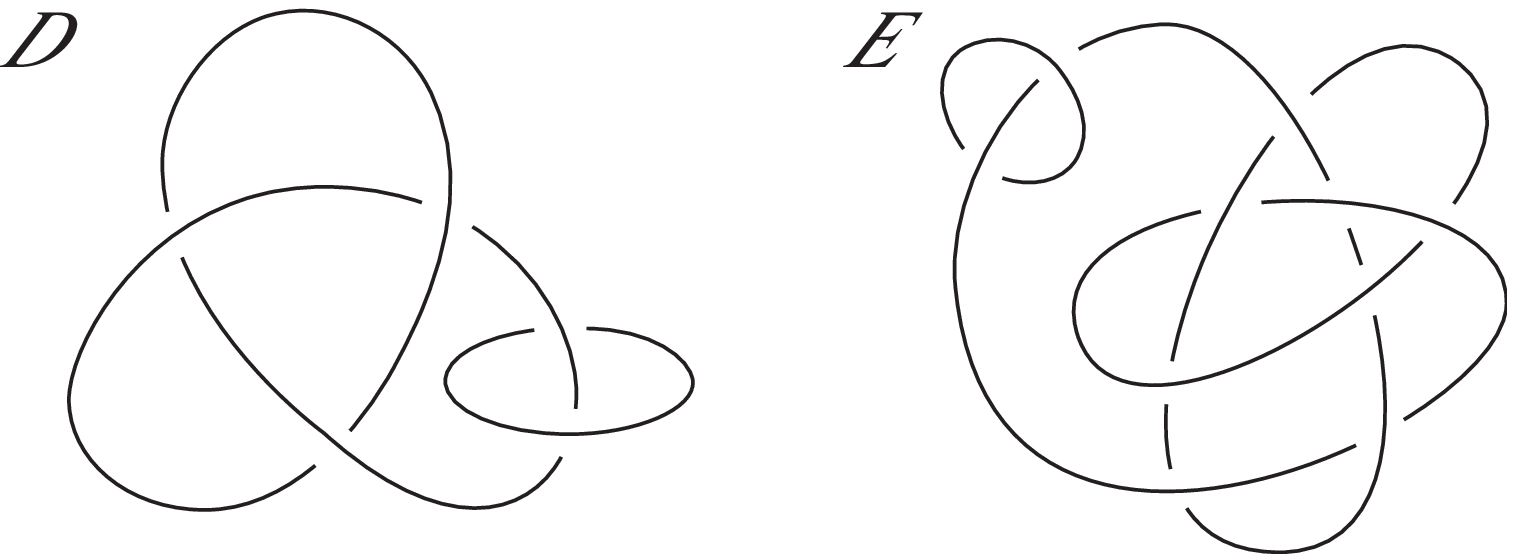}
\end{center}
\caption{}
\label{ex-fc}
\end{figure}
\end{eg}

\phantom{x}

\section{Relations of warping degree, unknotting number, and crossing number}

In this section, we enumerate several relations of the warping degree, the unknotting number or unlinking number, 
and the crossing number. 
Let $|D|$ be $D$ with orientation forgotten. 
We define the minimal warping degree of $D$ for all orientations as follows: 
$$d(|D|):=\min \{ d(D) | D: |D|\text{ with an orientation} \} .$$
Note that the minimal $d(|D|)$ for all diagrams $D$ of $L$ is equal to the ascending number $a(L)$ \cite{ozawa}: 
$$a(L)=\min \{ d(|D|) | D:\text{ a diagram of }L\} .$$
Let $E$ be a knot diagram, and $D$ a diagram of an $r$-component link. 
We review the relation of the unknotting number $u(E)$ (resp. the unlinking number $u(D)$) 
and the crossing number $c(E)$ (resp. $c(D)$) of $E$ (resp. $D$). 
The following inequalities are well-known \cite{nakanishi}: 
\begin{align}
u(E)\leq \frac{c(E)-1}{2},
\label{ucdk}
\end{align}

\begin{align}
u(D)\leq \frac{c(D)}{2}.
\label{ucdl}
\end{align}

\noindent Moreover, Taniyama mentioned the following conditions \cite{taniyama}: 

\phantom{x}

\noindent The necessary condition for the equality of (\ref{ucdk}) is that 
$E$ is a reduced alternating diagram of some $(2,p)$-torus knot, or $E$ is a diagram with $c(E)=1$.
The necessary condition for the equality of (\ref{ucdl}) is that 
every $D^i$ is a simple closed curve on $\mathbb{S}^2$ and every subdiagram $D^i\cup D^j$ is 
an alternating diagram. 

\phantom{x}

\noindent Hanaki and Kanadome characterized the link diagrams $D$ which satisfy $u(D)=(c(D)-1)/2$ as follows \cite{hanaki-kanadome}: 

\phantom{x}

\noindent Let $D=D^1\cup D^2\cup \dots \cup D^r$ be a diagram of an $r$-component link. 
Then we have 
$$u(D)=\frac{c(D)-1}{2}$$
if and only if exactly one of $D^1,D^2,\dots ,D^r$ is a reduced alternating diagram of a $(2,p)$-torus knot, 
the other components are simple closed curves on $\mathbb{S}^2$, and the non-self crossings of the subdiagram $D^i\cup D^j$ are 
all positive, all negative, or empty for each $i\neq j$. 
In addition, they showed that any minimal crossing diagram $D$ of a link $L$ with $u(L)=(c(L)-1)/2$ satisfies $u(D)=(c(D)-1)/2$. 

\phantom{x}

\noindent Abe and Higa study the knot diagrams $D$ which satisfy 
$$u(D)=\frac{c(D)-2}{2}.$$
Let $D$ be a knot diagram with $u(D)=(c(D)-2)/2$. 
They showed in \cite{abe-higa} that for any crossing point $p$ of $D$, 
one of the components of $D_p$ is a reduced alternating diagram of a $(2,p)$-torus knot 
and the other component of $D_p$ has no self-crossings, where $D_p$ is the diagram obtained from $D$ by smoothing at $p$. 
In addition, they showed that any minimal crossing diagram $D$ of a knot $K$ with $u(K)=(c(K)-2)/2$ satisfies the above condition. 

\phantom{x}

\noindent By adding to (\ref{ucdk}), we have the following corollary: 

\phantom{x}
\begin{cor}
For a knot diagram $E$, we have 
$$u(E)\leq d(|E|)\leq \frac{c(E)-1}{2}.$$
Further, if we have 
$$u(E)=d(|E|)=\frac{c(E)-1}{2},$$
then $E$ is a reduced alternating diagram of some $(2,p)$-torus knot, or $E$ is a diagram with $c(E)=1$.
\end{cor}
\phantom{x}

\noindent By adding to (\ref{ucdl}), we have the following corollary. 

\phantom{x}
\begin{cor}
\begin{description}
\item[(i)] For an $r$-component link diagram $D$, we have 
$$u(D)\leq d(|D|)\leq \frac{c(D)}{2}.$$
\item[(i\hspace{-1pt}i)] We have 
$$u(D)\leq d(|D|)=\frac{c(D)}{2}$$
if and only if every $D^i$ is a simple closed curve on $\mathbb{S}^2$ and the number of over-crossings of $D^i$ is equal to 
the number of under-crossings of $D^i$ in every subdiagram $D^i\cup D^j$ for each $i\neq j$. \\
\item[(i\hspace{-1pt}i\hspace{-1pt}i)] If we have 
$$u(D)=d(|D|)=\frac{c(D)}{2},$$
then every $D^i$ is a simple closed curve on $\mathbb{S}^2$ and for each pair $i$, $j$, the subdiagram $D^i\cup D^j$ is 
an alternating diagram. 
\end{description}

\phantom{x}
\begin{proof}
\begin{description}

\item[(i)] The equality $u(D)\leq d(|D|)$ holds because $u(D)\leq d(D)$ holds for every oriented diagram. 
We show that $d(|D|)\leq c(D)/2$. 
Let $D$ be an oriented diagram which satisfies 
$$d(D)=\sum ^r_{i=1} d(D^i)+ld(D)=d(|D|).$$
Then $D$ also satisfies 
\begin{align}
d(D^i)\leq \frac{c(D^i)}{2}
\label{101}
\end{align}
for every component $D^i$ because of the orientation of $D$. 
By Lemma \ref{2dc}, we have 
\begin{align}
ld(D)\leq \frac{lc(D)}{2}. 
\label{102}
\end{align}
Then we have 
$$\sum ^r_{i=1} d(D^i)+ld(D)\leq \sum ^r_{i=1} \frac{c(D^i)}{2}+\frac{lc(D)}{2}$$
by (\ref{101}) and (\ref{102}). Hence we obtain the inequality
$$d(|D|)\leq \frac{c(D)}{2}.$$\\

\item[(i\hspace{-1pt}i)] Suppose that the equality $d(|D|)=c(D)/2$ holds. 
Then the equalities
\begin{align}
d(D^i)=\frac{c(D^i)}{2}
\label{103}
\end{align}
and 
\begin{align}
ld(D)=\frac{lc(D)}{2}
\label{104}
\end{align}
hold by (\ref{101}) and (\ref{102}), where $D$ has an orientation such that $d(D)=d(|D|)$. 
The equality (\ref{103}) is equivalent to that $c(D^i)=0$ for every $D^i$. 
We prove this by an indirect proof. 
We assume that $c(D^i)>0$ for a component $D^i$. 
In this case, we have the inequality 
\begin{align}
d(D^i)+d(-D^i)+1\leq c(D^i)
\label{105}
\end{align}
by Theorem \ref{dk} since $D^i$ has a self-crossing. 
We also have 
\begin{align}
d(D^i)=d(-D^i)=\frac{c(D^i)}{2}
\label{106}
\end{align}
because $d(D^i)\leq d(-D^i)$ and (\ref{103}). 
By substituting (\ref{106}) for (\ref{105}), we have 
$$c(D^i)+1\leq c(D^i).$$
This implies that the assumption $c(D^i)>0$ is contradiction. 
Therefore every $D^i$ is a simple closed curve. 
The inequality (\ref{104}) is equivalent to that the number of over-crossings of $D^i$ is equal to 
the number of under-crossings of $D^i$ in every subdiagram $D^i\cup D^j$ for each $i\neq j$ by Lemma \ref{2dc}. 
On the other hand, suppose that every $D^i$ is a simple closed curve, and the number of over-crossings of $D^i$ is equal to 
the number of under-crossings of $D^i$ in every subdiagram $D^i\cup D^j$ for each $i\neq j$, then we have 
$$d(|D|)=ld(D)=\frac{lc(D)}{2}=\frac{c(D)}{2}.$$\\

\item[(i\hspace{-1pt}i\hspace{-1pt}i)] This holds by Corollary \ref{3cor}(i) and above Taniyama's condition. 

\end{description}
\end{proof}
\label{3cor}
\end{cor}

\phantom{x}

\noindent Let $K$ be a knot, and $L$ an $r$-component link. 
Let $u(K)$ be the unknotting number of $K$, and $u(L)$ be the unlinking number of $L$. 
The following inequalities are also well-known \cite{nakanishi}: 

\begin{align}
u(K)\leq \frac{c(K)-1}{2},
\label{uck}
\end{align}

\begin{align}
u(L)\leq \frac{c(L)}{2}.
\label{ucl}
\end{align}

\noindent The following conditions are mentioned by Taniyama \cite{taniyama}: 

\phantom{x}

\noindent The necessary condition for the equality of (\ref{uck}) is that $K$ is a ($2,p$)-torus knot ($p$:odd,$\neq \pm 1$). 
The necessary condition for the equality of (\ref{ucl}) is that 
$L$ has a diagram $D$ such that every $D^i$ is a simple closed curve on $\mathbb{S}^2$ 
and every subdiagram $D^i\cup D^j$ is an alternating diagram. 

\phantom{x}

\noindent By adding to (\ref{uck}) and (\ref{ek}), we have the following corollary: 

\phantom{x}
\begin{cor}
\begin{description}
\item[(i)] We have 
$$u(K)\leq \frac{e(K)}{2}\leq \frac{c(K)-1}{2}.$$
\item[(i\hspace{-1pt}i)] We have 
$$u(K)\leq \frac{e(K)}{2}=\frac{c(K)-1}{2}$$
if and only if $K$ is a prime alternating knot. \\
\item[(i\hspace{-1pt}i\hspace{-1pt}i)] If we have
$$u(K)=\frac{e(K)}{2}=\frac{c(K)-1}{2},$$
then $K$ is a ($2,p$)-torus knot ($p$:odd,$\neq \pm 1$). 
\end{description}
\end{cor}
\phantom{x}

\noindent By adding to (\ref{ucl}), we have the following corollary: 

\phantom{x}
\begin{cor}
For a diagram of an unoriented $r$-component link, we have 
$$u(L)\leq \frac{e(L)}{2}\leq \frac{c(L)}{2}.$$
Further, if the equality $u(L)=e(L)/2=c(L)/2$ holds, then 
$L$ has a diagram $D=D^1\cup D^2\cup \dots \cup D^r$ such that every $D^i$ is a simple closed curve on $\mathbb{S}^2$ 
and for each pair $i,j$, the subdiagram $D^i\cup D^j$ is an alternating diagram. 

\phantom{x}
\begin{proof}
We prove the inequality $u(L)\leq e(L)/2$. 
Let $D$ be a minimal crossing diagram of $L$ which satisfies $e(L)=d(D)+d(-D)$. 
Then we obtain 
$$e(L)=d(D)+d(-D)\geq 2u(D)\geq 2u(L).$$
The condition for the equality is due to above Taniyama's condition.
\end{proof}
\end{cor}
\phantom{x}

\section{Splitting number}

In this section, we define the splitting number and enumerate relations of the warping degree and the complete splitting number. 
The \textit{splitting number} (resp. \textit{complete splitting number}) of $D$, 
denoted by $Split(D)$ (resp. $split(D)$), is the smallest number of crossing changes 
which are needed to obtain a diagram of a splittable (resp. completely splittable) link from $D$. 
The splitting number of a link which is the minimal $Split(D)$ for all diagrams $D$ is defined by Adams \cite{adams}. 
The \textit{linking splitting number} (resp. \textit{linking complete splitting number}) of $D$, 
denoted by $lSplit(D)$ (resp. $lsplit(D)$), is the smallest number of non-self-crossing changes 
which are needed to obtain a diagram of a splittable (resp. completely splittable) link from $D$. 
We have the following propositions: 

\phantom{x}
\begin{prop}
\begin{description}
\item[(i)] We have 
$$split(D)\leq d(|D|).$$\\

\item[(i\hspace{-1pt}i)] We have 
$$split(D)\leq lsplit(D)\leq ld(D)\leq \frac{lc(D)}{2}\leq \frac{c(D)}{2}.$$
\end{description}
\label{lsplit}
\end{prop}

\phantom{x}

\noindent We give examples of Proposition \ref{lsplit}.

\begin{eg}
The diagram $D$ in Figure \ref{ex-sd} has $split(D)=2<d(|D|)=3$. 
The diagram $E$ in Figure \ref{ex-sd} has $split(E)=d(|E|)=3$. 
\begin{figure}[h]
\begin{center}
\includegraphics[width=70mm]{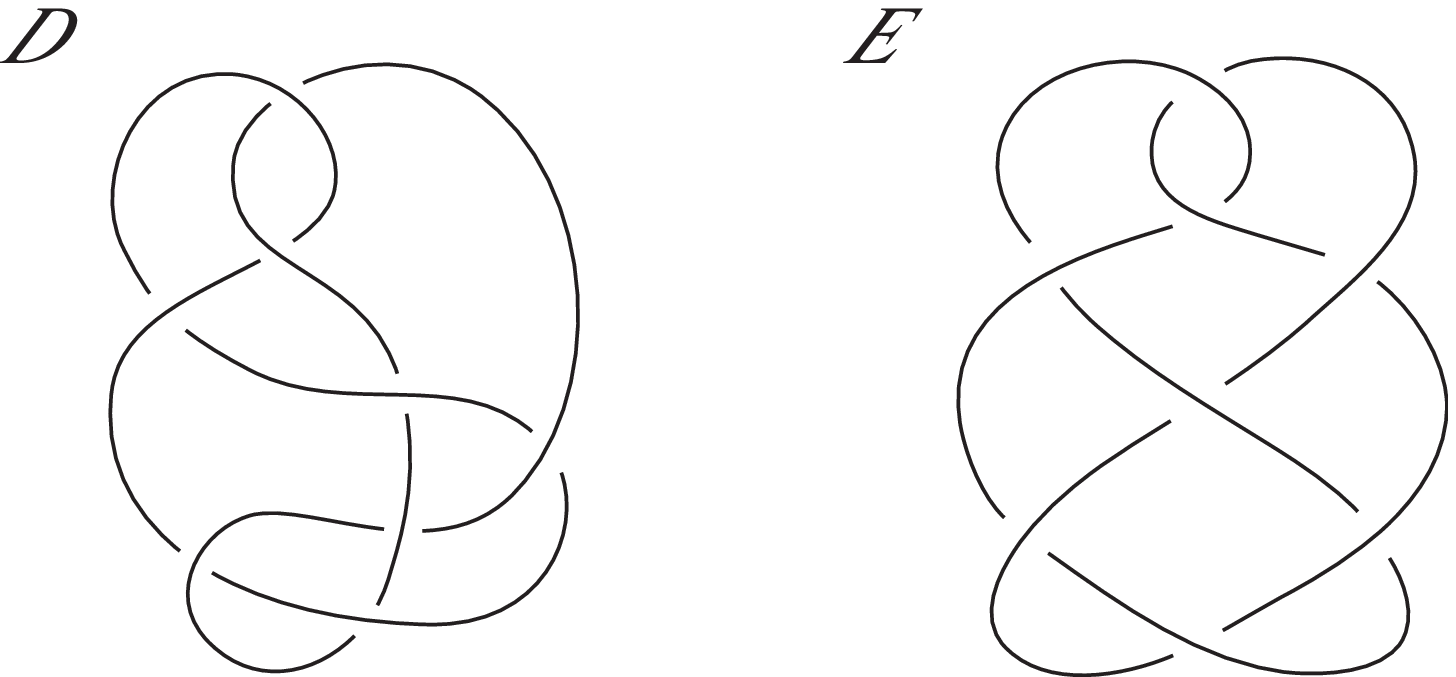}
\end{center}
\caption{}
\label{ex-sd}
\end{figure}
\end{eg}

\phantom{x}
\begin{eg}
The diagram $D$ in Figure \ref{ex-splitl} has $split(D)=1<lsplit(D)=2$. 
The diagram $E$ in Figure \ref{ex-splitl} has $split(E)=lsplit(E)=2$. 
\begin{figure}[h]
\begin{center}
\includegraphics[width=65mm]{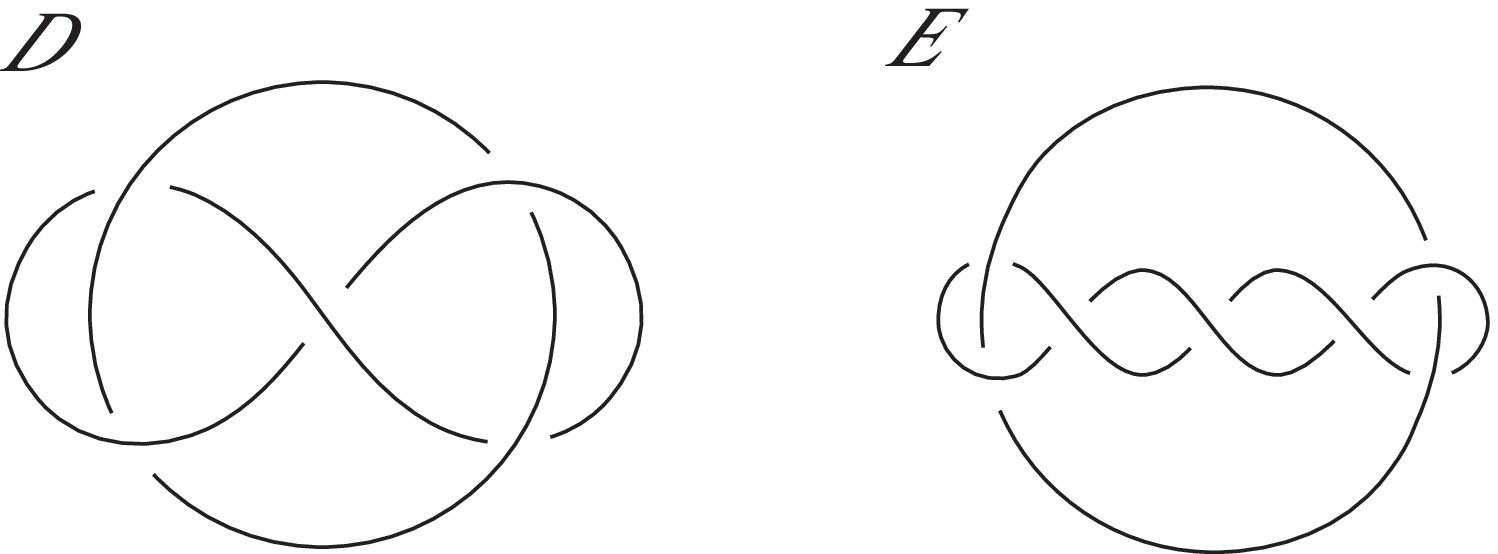}
\end{center}
\caption{}
\label{ex-splitl}
\end{figure}
\end{eg}

\phantom{x}
\begin{eg}
The diagram $D$ in Figure \ref{ex-lsplitld} has $lsplit(D)=3<ld(D)=5$. 
The diagram $E$ in Figure \ref{ex-lsplitld} has $lsplit(E)=ld(E)=5$. 
\begin{figure}[h]
\begin{center}
\includegraphics[width=80mm]{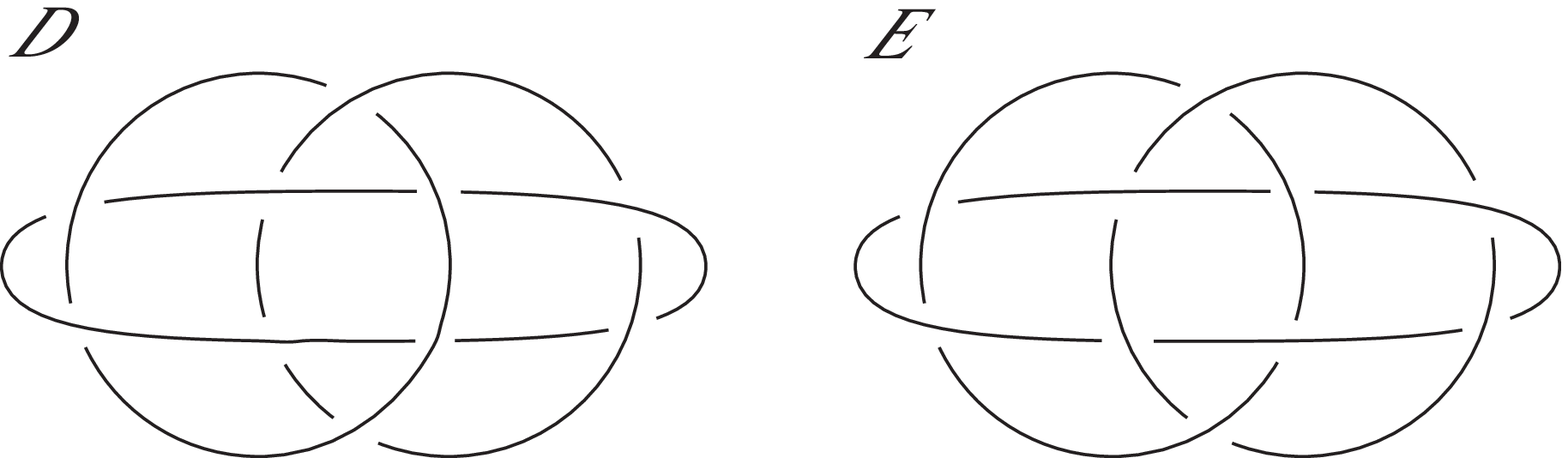}
\end{center}
\caption{}
\label{ex-lsplitld}
\end{figure}
\end{eg}
\phantom{x}

\noindent We raise the following question: 

\phantom{x}
\begin{q}
When does the equality 
$$split(D)=d(|D|),$$
$$split(D)=lsplit(D)$$
or
$$lsplit(D)=ld(D)$$
hold?
\end{q}
\phantom{x}

\section{Calculation of warping degree}

In this section, we show methods for calculating the warping degree and linking warping degree by using matrices. 
First, we give a method for calculating the warping degree $d(D)$ of an oriented knot diagram $D$. 
Let $a$ be a base point of $D$. 
We can obtain the warping degree $d(D_a)$ of $D_a$ by counting the warping crossing points easily. 
Let $[D_a]$ be a sequence of some "$o$" and "$u$", which is obtained as follows. 
When we go along the oriented diagram $D$ from $a$, 
we write down "$o$" (resp. "$u$") if we reach a crossing point as an over-crossing (resp. under-crossing) in numerical order. 
We next perform normalization to $[D_a]$, by deleting the subsequence "$ou$" repeatedly, 
to obtain the normalized sequence $\lfloor D_a\rfloor $. 
Then we have 
$$d(D)=d(D_a)-\frac{1}{2}\sharp \lfloor D_a\rfloor ,$$
where $\sharp \lfloor D_a\rfloor $ denotes the number of entries in $\lfloor D_a\rfloor $. 
Thus, we obtain the warping degree $d(D)$ of $D$. 
In the following example, we find the warping degree of a knot diagram by using the above algorithm. 

\begin{eg}
For the oriented knot diagram $D$ and the base point $a$ in Figure \ref{948}, 
we have $d(D_a)=4$ and $[D_a]=[oouuouuouuouoouoou]$. 
By normalizing $[D_a]$, we obtain $\lfloor D_a\rfloor =[uuoo]$. 
Hence we find the warping degree of $D$ as follows: 
$$d(D)=4-\frac{1}{2}\times 4=2.$$
\begin{figure}[h]
\begin{center}
\includegraphics[width=30mm]{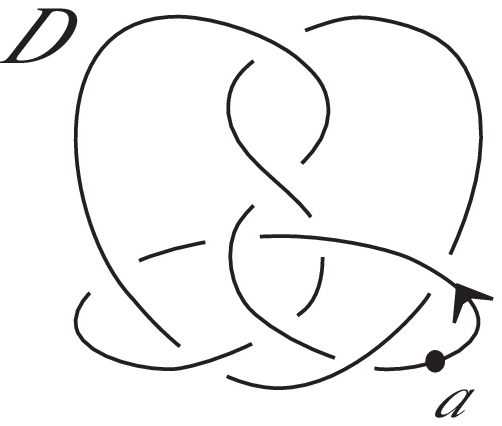}
\end{center}
\caption{}
\label{948}
\end{figure}
\end{eg}

\noindent For some types of knot diagram, this algorithm is useful in formulating the warping degree or looking into its properties. 
We enumerate the properties of an oriented diagram of a pretzel knot of odd type in the following example: 

\begin{eg}
Let $D=P(\varepsilon _1n_1, \varepsilon _2n_2, \dots ,\varepsilon _mn_m)$ be an oriented pretzel knot diagram of odd type 
($\varepsilon _i\in{+1,-1}, n_i,m$: odd$>0$), where the orientation is given as shown in Figure \ref{pretzel}. 
We take base points $a$, $b$ in Figure \ref{pretzel}. 
Then we have 
$$d(D_a)=d(-D_b)=\frac{c(D)}{2}+\sum _i \frac{(-1)^{i+1}\varepsilon _i}{2}$$
and 
$$\sharp \lfloor D_a\rfloor =\sharp \lfloor -D_b\rfloor .$$
Hence we have $d(D)=d(-D)$ in this case. 
In particular, if $D$ is alternating i.e. $\varepsilon _1=\varepsilon _2=\dots =\varepsilon _m=\pm 1$, then we have that 
$$d(D)=\frac{c(D)}{2}-\frac{1}{2}.$$
\begin{figure}[h]
\begin{center}
\includegraphics[width=40mm]{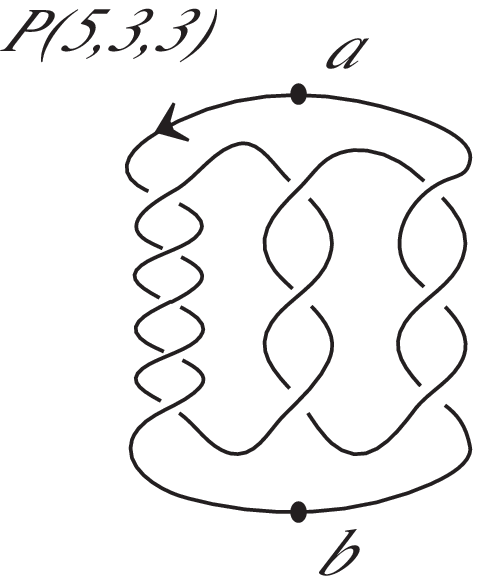}
\end{center}
\caption{}
\label{pretzel}
\end{figure}
\end{eg}

\noindent We next explore how to calculate the linking warping degree $ld(D)$ by using matrices. 
For a link diagram $D$ and a base point sequence ${\bf a}$ of $D$, 
we define an $r$-square matrix $M(D_{\bf a})=(m_{i j})$ by the following rule: 

\begin{itemize}
\item For $i\neq j$, $m_{i j}$ is the number of crossings of $D^i$ and $D^j$ which are under-crossings of $D^i$. 
\item For $i=j$, $m_{i j}=d(D^i)$. 
\end{itemize}

\noindent We show an example. 

\phantom{x}
\begin{eg}
For $D_{\bf a}$ and $D_{\bf b}$ in Figure \ref{matrix}, we have

\begin{eqnarray*}
M(D_{\bf a})=
\left(
\begin{array}{ccc}
0&1&0\\
1&0&0\\
2&2&0\\
\end{array}
\right) ,\ M(D_{\bf b})=
\left(
\begin{array}{ccc}
0&2&2\\
0&0&1\\
0&1&0\\
\end{array}
\right) .
\end{eqnarray*}

\begin{figure}
\begin{center}
\includegraphics[width=70mm]{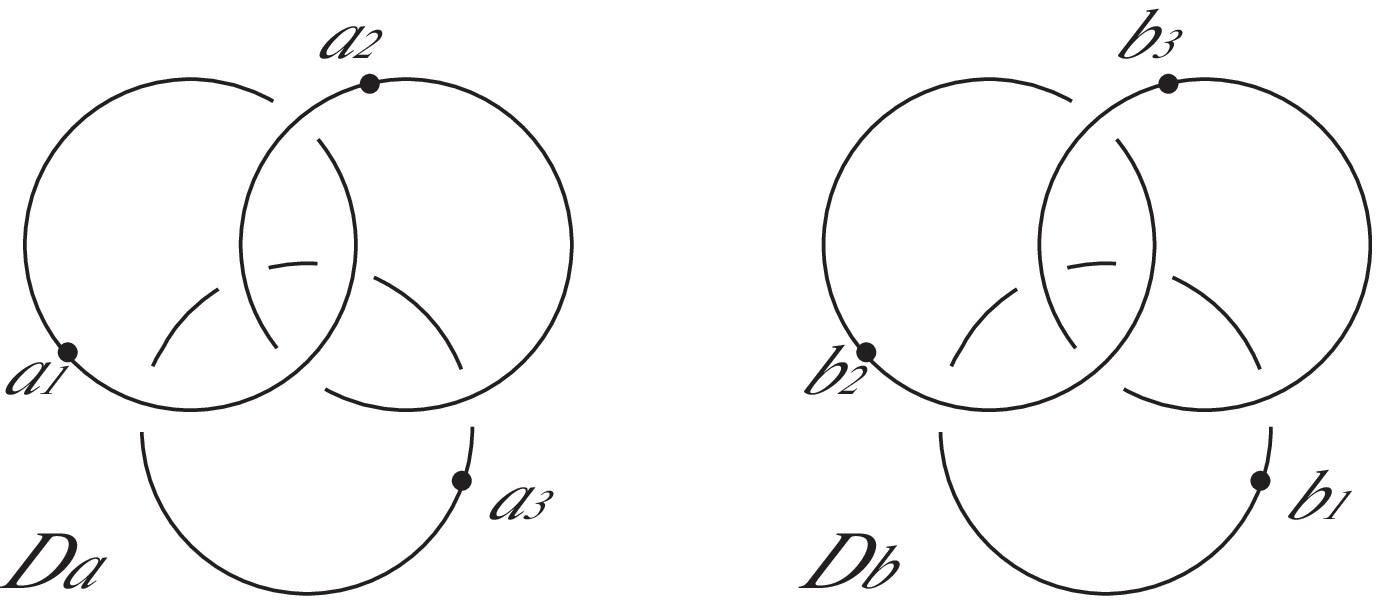}
\end{center}
\caption{}
\label{matrix}
\end{figure}

\end{eg}
\phantom{x}

\noindent We note that $ld(D_{\bf a})$ is obtained by summing the upper triangular entries of $M(D_{\bf a})$, that is 
$$ld(D_{\bf a})=\sum _{i<j}m_{i j},$$
and we notice that 
$$d(D_{\bf a})=\sum _{i\leq j}m_{i j},$$
where $m_{i j}$ is an entry of $M(D_{\bf a})$ ($i,j=1,2,\dots ,r$). 
For the base point sequence ${\bf a}'=(a_1,a_2,\dots ,a_{k+1},a_k,\dots ,a_r)$ which is obtained from a base point sequence ${\bf a}$ 
by exchanging $a_k$ and $a_{k+1}$ ($k=1,2,\dots ,r-1$), the matrix $M(D_{{\bf a}'})$ is obtained as follows: 
$$M(D_{{\bf a}'})=P_kM(D_{\bf a})P_k^{-1},$$
where
\begin{eqnarray*}
P_k=
\left(
\begin{array}{cccccccc}
1&&&&&&&\\
&\ddots &&&&&&\\
&&&0&1&&&\\
&&&1&0&&&\\
&&&&&&\ddots &\\
&&&&&&&1\\
\end{array}
\right) 
; m_{i j}=\left\{
\begin{array}{l}
1\text{ for }(i,j)=(k,k+1),(k+1,k)\\
     \hspace{3mm}\text{and}(i,j)=(i,i) (i\neq k,k+1),\\
0\text{ otherwise. }
\end{array}
\right.
\end{eqnarray*}
\phantom{x}

\noindent With respect to the linking warping degree, we have 
$$ld(D_{{\bf a}'})=ld(D_{\bf a})-m_{k k+1}+m_{k+1 k},$$
where $m_{k k+1}, m_{k+1 k}$ are entries of $M(D_{\bf a})$. 
To enumerate the permutation of the order of ${\bf a}=(a_1,a_2,\dots ,a_r)$, we consider a matrix 
$Q=P^{r-1}P^{r-2}\dots P^2P^1$, where $P^n$ denotes $P_nP_{n+1}\dots P_{k_n}$ ($n\leq k_n \leq r-1$) or the identity matrix $E_r$. 
Since $Q$ depends on the choices of $k_n$ ($n=1,2,\dots ,r-1$), 
we also denote $Q$ by $Q_{\mathbf{k}}$, 
where $\mathbf{k}=(k_1,k_2,\dots ,k_{r-1})$ ($n\leq k_n \leq r$) and we regard $P^n=E_r$ in the case $k_n=r$. 
Hence we obtain the following formula: 
$$ld(D)=\min _{\mathbf{k}} \{ \sum _{i<j}m_{i j} | m_{i j}:\text{ an entry of }Q_{\mathbf{k}}M(D_{\bf a})Q_{\mathbf{k}}^{-1}\}.$$
Thus, we obtain the warping degree of an oriented link diagram by summing the warping degrees $d(D^i)$ ($i=1,2,\dots ,r$) 
and the linking warping degree $ld(D)$.



\maketitle

\end{document}